# Non Additive Geometry

Shai Haran

shaiharan@gmail.com

# Contents

# Preface

The usual dictionary between geometry and commutative algebra is not appropriate for Arithmetic geometry because addition is a singular operation at the "Real prime". We replace Rings, with addition and multiplication, by Props (=strict symmetric monoidal category generated by one object), or by Bioperad (=two closed symmetric operads acting on each other): to a ring we associate the prop of all matrices over it, with matrix multiplication and block direct sums as the basic operations, or the bioperad consisting of all raw and column vectors over it. We define the "commutative" props and bioperads, and using them we develop a generalized algebraic geometry, following Grothendieck footsteps closely. This new geometry is appropriate for Arithmetic (and potentially also for Physics).

# CHAPTER 0

# Introduction

We, the human beings, are walking on this earth doing most of the time addition and subtraction in our brains - "How much will I profit? earn? pay? How much it will cost? What will be left?" . There's no wonder that addition and subtraction dominate all our mathematics. It is always an Abelian group that is at the basis of any mathematical structures (or an additive functor between Abelian categories, and homological algebra that produces all the deeper theorems). However, there is a world before addition, before we apply the functor of "the free Abelian group" and start doing addition and subtraction - it is the world of homotopy. It is the long exact sequences of (co)fibrations that give us all the deeper theorems we need.

The encoding of geometry into an Abelian group begins with Gelfand theorem on the equivalence of the category of compact Hausdorff topological spaces and continuous maps, with the opposite of the category of commutative $C^*$- algebra and $*$-homomorphism. To a topological space $X$ one associate the commutative $C^*$- algebra $C(X) = \{f : X \to \mathbb{C} \text{ continuous}\}$ using the addition and multiplication (and conjugation and norm) in $\mathbb{C}$ to get the structure of a $C^*$-algebra on $C(X)$. To a commutative $C^*$-algebra A one associate the space spec($A$) of maximal ideals of $A$, which is compact Hausdorff with respect to the minimal topology making all $f \in A$ continuous.

Alexander Grothendieck (Shapiro) came to algebraic geometry from $C^*$-algebras, and in front of him were the Weil conjectures which are all "above $\mathbb{Z}$". He therefore chose **Commutative Rings** as the basis for the language of algebraic geometry. If he was more interested in Arithmetic he would have had to choose a more general notion than a commutative ring. That commutative rings are not sufficient for Arithmetic geometry is clear by considering Weil's "Rosetta stone", that's the one dimensional structures we have:

(1*)   Number Fields: $K/\mathbb{Q}$ finite.

(2*)   Function Fields in characteristic $p$: $K/\mathbb{F}_q(z)$ finite.

(3*)   Function Fields in characteristic 0: $K/\mathbb{C}(z)$ finite.

The Function Fields $K$ corresponds one to one with smooth projective curves $X_K$, and the points of $X_K$ are in bijection with the valuation sub rings of $K$. A choice of the transcendental basis "$z$" corresponds to a choice of a dominant mapping $z : X_K \twoheadrightarrow \mathbb{P}^1$. At the basis of the Rosetta stone we have the rational numbers/functions, with their real prime/point at infinity:

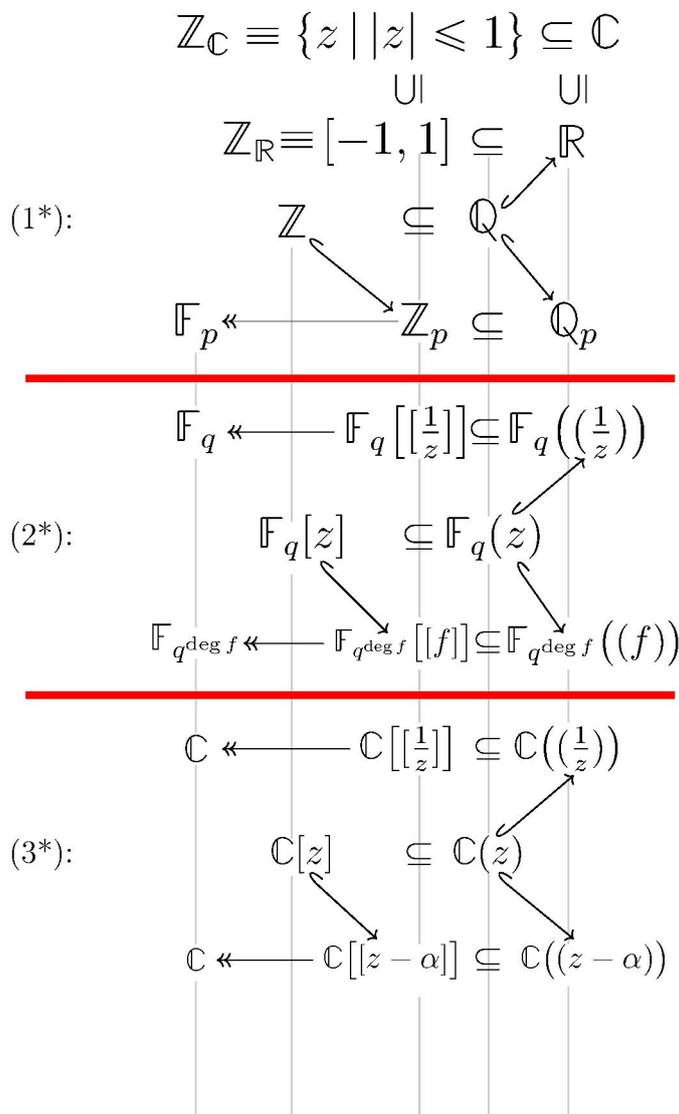

We see that the "Real integers" $\mathbb{Z}_{\mathbb{R}} = [-1, 1] \subseteq \mathbb{R}$, and similarly the "Complex integers" $\mathbb{Z}_{\mathbb{C}} = \{z \in \mathbb{C}, |z| \leq 1\}$, are not closed under addition, and do not form a sub ring. As Grothendieck knew well: the global theorems exist only in projective geometry, not in affine geometry. The most basic example being that a rational function without poles anywhere, including the point at infinity, must be a constant. If our language does not "see" all of the points there will be no global theorems. We must therefore understand these real and complex "integers" as some kind of a "generalized ring".

Another problem is the mysterious "Arithmetic surface": the integers $\mathbb{Z}$ are the initial object of the category of rings, so the categorical sum of $\mathbb{Z}$ with itself reduces to the diagonal

$$(0.1) \qquad \mathbb{Z} \otimes \mathbb{Z} = \mathbb{Z},$$

unlike the function field analog

$$(0.2) \qquad F[z] \otimes_F F[z] = F[z_1, z_2].$$

Thus if we insist on rings, the **Arithmetical Surface**

$$(0.3) \qquad \operatorname{spec}(\mathbb{Z}) \coprod \operatorname{spec}(\mathbb{Z}) = \operatorname{spec}(\mathbb{Z})$$

reduces to its diagonal.

The original Rosetta Stone had three languages talking word for word about one and the same reality. Comparing the well understood ancient Greek and Demotic languages with the mysterious Hieroglyphic, the ancient Egyptian language, led to the decipherment of the Hieroglyphics, and with it the understanding of many texts from the mummies and the pyramids. While in mathematics we have one language that is talking about three "one dimensional" but different realities (1*), (2*), and (3*). Moreover, this language of commutative rings, the language of **Addition and Multiplication**, is not the right language for Arithmetic geometry, because addition is singular at the real prime, and the Arithmetic Surface disappears.

A similar situation occurs in physics: normalizing the speed of light to one, $c = 1$, the interval of speeds $[-1, 1]$ is not closed under addition. Einstein's solution was to replace $z_1 + z_2$ by

$$(0.4) \qquad z_1(+)z_2 := (z_1 + z_2)/(1 + z_1 \bullet z_2)$$

That is we use the fractional linear transformation $z \mapsto (1+z)/(1-z)$ that identifies the interval $[-1, 1]$ with the positive reals $[0, \infty]$, and we carry multiplication of positive real numbers $(0, \infty)$ to the interval $(-1, 1)$ to obtain Einstein's addition $(+)$. When we do it with complex numbers the unit disc $\mathbb{Z}_{\mathbb{C}} = \{|z| \leq 1\}$ get identified with the right half plane $\{Re(z) \geq 0\}$ which is not closed with respect to multiplication, so $\mathbb{Z}_{\mathbb{C}}$ is not closed under $(+)$. (Curiously, we do have that $\mathbb{Z}_{\mathbb{C}}$ is closed with respect to the operations $z_1(+)'z_2 := (z_1 + \bar{z}_2)/(1 + z_1 z_2)$, or $(z_1 + z_2)/(1+z_1\bar{z}_2)$ but these operations are not commutative and are not associative). It is perhaps already here that relativity and quantum mechanics cannot coexist. Quantum mechanics started with the observation by Heisenberg that the level structure of the energy of microscopic particles are controlled by finite matrices.

Our basic idea is to replace the ring $A$ by the collection of all matrices over $A = \{A_{n,m}\}$, $n, m \geq 0$ with the operations of **Matrix Multiplication and of Block Direct Sum of Matrices**. The algebraic structure of this collection of all matrices is that of a "Prop": a strict symmetric monoidal category generated by one object.

Observing that

$$(0.5) \qquad (\mathbb{Z}_p)_{n,m} = \{a \in (\mathbb{Q}_p)_{n,m} \mid a \circ \mathbb{Z}_p^m \subseteq \mathbb{Z}_p^n\}$$

we define

$$(0.6) \qquad (\mathbb{Z}_{\mathbb{R}})_{n,m} = \{a \in (\mathbb{R})_{n,m} \mid a \circ \mathbb{Z}_{\mathbb{R}}^m \subseteq \mathbb{Z}_{\mathbb{R}}^n\}$$

with

$$(0.7) \qquad \mathbb{Z}_{\mathbb{R}}^m = \{(z_1, \ldots, z_m) \in \mathbb{R}^m, |z_1|^2 + \cdots + |z_m|^2 \leq 1\}$$

the unit $\ell_2$ ball. Similarly the complex integers are defined by

$$(0.8) \qquad (\mathbb{Z}_{\mathbb{C}})_{n,m} = \{a \in \mathbb{C}_{n,m} \mid a \circ \mathbb{Z}_{\mathbb{C}}^m \subseteq \mathbb{Z}_{\mathbb{C}}^n\}$$

with

$$(0.9) \qquad \mathbb{Z}_\mathbb{C}^m = \{(z_1, \ldots, z_m) \in \mathbb{C}^m, |z_1|^2 + \ldots + |z_m|^2 \leq 1\}$$

the unit $\ell_2$ complex ball. The sub props

$$(0.10) \qquad Z_\mathbb{R} \subseteq \mathbb{R}, \qquad Z_\mathbb{C} \subseteq \mathbb{C}$$

are closed under the operations of matrix multiplication and block direct sum of matrices.

We can define the "residue field at the real prime" $\mathbb{F}_\mathbb{R}$, and similarly the "residue field at the complex prime" $\mathbb{F}_\mathbb{C}$, as follows. For $a \in (\mathbb{Z}_\mathbb{R})_{n,m}$ we have $a^t \in (\mathbb{Z}_\mathbb{R})_{m,n}$, and the self-adjoint operators $a^t \circ a \in (\mathbb{Z}_\mathbb{R})_{m,m}$, $a \circ a^t \in (\mathbb{Z}_\mathbb{R})_{n,n}$. The spectral decomposition of these operators give

$$(0.11) \qquad \mathbb{R}^n = \ker(a^t) \oplus \bigoplus_{i=1}^k W(\lambda_i) \qquad \mathbb{R}^m = \ker(a) \oplus \bigoplus_{i=1}^k V(\lambda_i)$$

with $0 < \lambda_1 < \lambda_2 < \cdots < \lambda_k \leq 1$, and $a$ induces isomorphism
$$(0.12)$$
$$W(\lambda_i) = \ker\left(a \circ a^t - \lambda_i \circ \operatorname{Id}_n\right) \xleftarrow[\sim]{a} V(\lambda_i) = \ker\left(a^t \circ a - \lambda_i \cdot \operatorname{Id}_m\right)$$

In particular, $a$ induces the partial isometry (possibly empty!)

$$(0.13) \qquad \widehat{a} := \{W(1) \xleftarrow[\sim]{a} V(1)\}$$

Let $\mathbb{F}_\mathbb{R}$ denote the collection of all such partial isometries

$$(F_\mathbb{R})_{n,m} := \{\mathcal{W}(a) \xleftarrow[\sim]{a} \mathcal{V}(a) \text{ linear isometry}, \mathcal{W}(a) \subseteq \mathbb{R}^n, \mathcal{V}(a) \subseteq \mathbb{R}^m\}.$$

The composition in $\mathbb{F}_\mathbb{R}$ is given by
$$(0.14)$$
$$\left(\mathcal{W}(a) \xleftarrow[\sim]{a} \mathcal{V}(a)\right) \circ \left(\mathcal{W}(b) \xleftarrow[\sim]{b} \mathcal{V}(b)\right) := \left(a(\mathcal{V}(a) \cap \mathcal{W}(b)) \xleftarrow[\sim]{a \circ b} b^{-1}(\mathcal{V}(a) \cap \mathcal{W}(b))\right)$$

The direct sum in $\mathbb{F}_\mathbb{R}$ is the usual sum

(0.15)
$$\left(\mathcal{W}(a_1)) \xleftarrow[\sim]{a_1} \mathcal{V}(a_1)\right) \oplus \left(\mathcal{W}(a_2) \xleftarrow[\sim]{a_2} \mathcal{V}(a_2)\right) :=$$
$$\left(\mathcal{W}(a_1) \oplus \mathcal{W}(a_2) \xleftarrow[\sim]{a_1 \oplus a_2} \mathcal{V}(a_1) \oplus \mathcal{V}(a_2)\right)$$

With these operations $\mathbb{F}_\mathbb{R}$ forms a prop, and we have the surjection

(0.16)
$$\mathbb{Z}_\mathbb{R} \twoheadrightarrow \mathbb{F}_\mathbb{R}$$
$$a \mapsto \widehat{a} := \left(W(1) \xleftarrow[\sim]{a} V(1)\right)$$

In exactly the same way (using $\overline{a}^t$ instead of $a^t$) we have the surjection of props

(0.17)
$$\mathbb{Z}_\mathbb{C} \twoheadrightarrow \mathbb{F}_\mathbb{C}$$

With these definitions, we see that the "mathematical Rosetta stone" can be repaired, and in the language of props, the Arithmetic reality (1*) becomes word for word analogue to the function fields realities (2*) and (3*):

(0.18)
$$\begin{array}{ccccc}
\mathbb{C} & \supseteq & \mathbb{Z}_\mathbb{C} & \twoheadrightarrow & \mathbb{F}_\mathbb{C} \\
\cup\mathsf{I} & & \cup\mathsf{I} & & \cup\mathsf{I} \\
\mathbb{R} & \supseteq & \mathbb{Z}_\mathbb{R} & \twoheadrightarrow & \mathbb{F}_\mathbb{R} \\
& & \cup\mathsf{I} & & \\
& & \mathbb{Z} & \subseteq & \mathbb{Q} \\
& & \cap & & \cap \\
\mathbb{F}_p & \twoheadleftarrow & \mathbb{Z}_p & \subseteq & \mathbb{Q}_p
\end{array}$$

Notice that the prop associated to any ring A will contain the sub prop $\mathbb{F}$ with

(0.19)  $\mathbb{F}_{n,m} = $ | all the n by m matrices with entries 0 or 1 such that in any row and in any column there will be at most one nonzero entry. |

This common sub prop of all rings is the **"Field with One element"** [**Sou08**]. It will be the initial object of our props, and the **"absolute-point"** spec($\mathbb{F}$) will be the last object of our geometry.

There is no such thing as non commutative geometry ! Although non commutative rings (and similarly non-commutative props) have a rich theory of cohomology and $K$-theory generalizing the theory of commutative rings [**Con94**], there is no spectrum associated with a non commutative algebra, there are no spaces, localization and sheaves. We will have to understand what it means for the prop A when the underlying ring A is commutative. There will be two notions of commutativity for props giving rise to the full embeddings beginning with the category of commutative rings $\mathcal{CRing}$,

(0.20) $$\mathcal{CRing} \subseteq \mathcal{C}_T\mathcal{Prop} \subseteq \mathcal{CProp} \subseteq \mathcal{Prop}.$$

These categories are all complete and co-complete. In particular we have the Arithmetical Surface , the categorical sum of $\mathbb{Z}$ with itself in the category $\mathcal{CProp}$

(0.21) $$\begin{array}{c} \mathbb{Z} \otimes_{\mathbb{F}} \mathbb{Z} \twoheadrightarrow \mathbb{Z} \\ \searrow \quad \nearrow \\ (\mathbb{Z} \otimes_{\mathbb{F}} \mathbb{Z})^T \end{array}$$

Unfortunately, the same sum in the category $\mathcal{C}_T\mathcal{Prop}$, the total commutative quotient of $\mathbb{Z} \otimes_{\mathbb{F}} \mathbb{Z}$, reduces to $\mathbb{Z}$. This is why we prefer to work with the weaker notion of commutativity.

To a commutative prop $A$ we will associate functorially a compact topological space spec($A$). This space depends on much less information than what $A$ gives. To construct spec($A$) we only need the information given by the collection of "raw and column vectors in $A$", we call this $UA$, so

(0.22) $\quad (UA)^-(n) := A_{1,n} \quad , \quad (UA)^+(n) = A_{n,1}$

The collection $UA = (P^-(n), P^+(n))$ forms the structure we call Bio, short for Bi-operad, consisting of two (closed symmetric) operads $P^-$ and $P^+$ acting on each other in a consistent manner. We obtain a commutative diagram with rows full embeddings and their left adjoints

(0.23)
$$\begin{array}{c}
C_T \mathcal{P}\!\mathit{rop} \leftrightarrows C\mathcal{P}\!\mathit{rop} \leftrightarrows \mathcal{P}\!\mathit{rop} \\
C\mathcal{R}\mathit{ing} \quad \mathcal{F} \uparrow \downarrow U \quad \mathcal{F} \uparrow \downarrow U \quad \mathcal{F} \uparrow \downarrow U \\
C_T \mathcal{B}\mathit{io} \leftrightarrows C\mathcal{B}\mathit{io} \leftrightarrows \mathcal{B}\mathit{io}
\end{array}$$

Here $C\mathcal{B}\mathit{io}$ are the bios $P = (P^-, P^+)$ where all operations $p \in P^-(n)$, and all co-operations $q \in P^+(m)$, commute in the sense of bi-algebras,

(0.24) $\quad q \circ p = \underbrace{(p, \cdots, p)}_{m} \circ \sigma_{m,n} \circ \underbrace{(q, \cdots, q)}_{n}.$

Thus we identify

(0.25)

```
        1              1
         \              \
          \ :            \ :
     i ----•q         p•---- j
          / :            / :
         /              /
        m              n
```

with

(0.26) 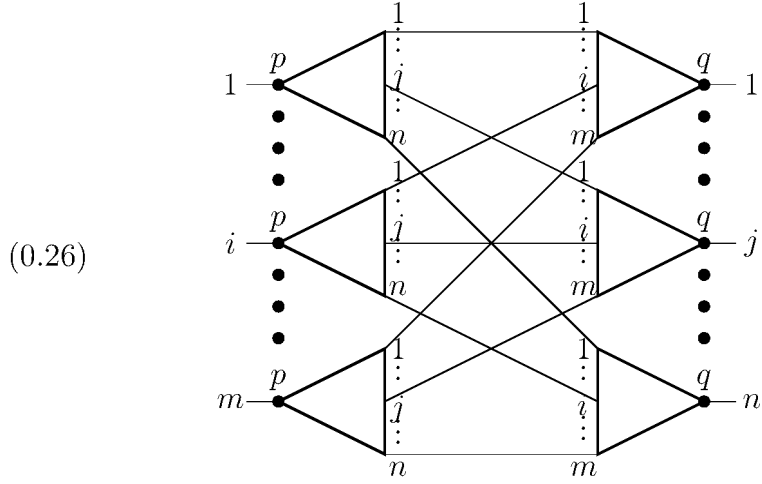

$$(0.27) \quad (i-1)\cdot n + j = \sigma_{m,n}\bigl((j-1)m+i\bigr), \quad 1 \leqslant i \leqslant m, \ 1 \leqslant j \leqslant n,$$

the permutation $\sigma_{m,n} \in S_{m\bullet n}$ is given pictorially in the diagram (0.26). The totally commutative bios $\mathcal{C}_T\mathcal{B}\mathit{io}$ are the commutative bios $P = (P^-, P^+)$ where the operads $P^-$ and $P^+$ satisfy further the original Boardman-Vogt commutativity:

(0.28)
$$\text{for } p \in P^-(m),\ p' \in P^-(n)\ :\ p \circ (\underbrace{p',\cdots,p'}_{m}) = p' \circ (\underbrace{p,\cdots,p}_{n}) \circ \sigma_{n,m}$$
$$\text{for } q \in P^+(m),\ q' \in P^+(n)\ :\ (\underbrace{q',\cdots,q'}_{m}) \circ q = \sigma_{m,n} \circ (\underbrace{q,\cdots,q}_{n}) \circ q'$$

For a commutative bio $P = (P^-, P^+)$, (resp. prop $A$), the monoid $P^-(1) \cong P^+(1)\ (= A_{1,1})$ is commutative and central, so that we can localize with respect to multiplicative subsets $S \subseteq P(1)\ (= A_{1,1})$. We obtain a sheaf $\mathcal{O}_P$ (resp. $\mathcal{O}_A$) of commutative bios (resp. props) with local stalks over spec($P$) (resp. spec(U$A$)). To a map of bios (props) $\varphi$ corresponds a map of locally bio (resp. prop) spaces in the opposite direction. We obtain the diagram of full-embedding

(0.29)
$$\begin{array}{ccc} (\mathcal{CRing})^{op} & \longleftrightarrow \mathcal{S}ch & \longleftrightarrow \mathcal{CRing}/\mathcal{T}op \\ \updownarrow & \updownarrow & \updownarrow \\ (\mathcal{CBis})^{op} & \longleftrightarrow \mathcal{BSch} & \longleftrightarrow \mathcal{CBis}/\mathcal{T}op \\ \updownarrow & \updownarrow & \updownarrow \\ (\mathcal{CProp})^{op} & \longleftrightarrow \mathcal{PSch} & \longleftrightarrow \mathcal{CProp}/\mathcal{T}op \end{array}$$

Here the column on the right consists of pairs $(X, \mathcal{O}_X)$ of a topological space $X$ together with a sheaf $\mathcal{O}_X$ of $\mathcal{CRing}/\mathcal{CBis}/\mathcal{CProp}$ with local stalks at each point $x \in X$

(0.30) $$\mathcal{O}_{X,x} := \varinjlim_{x \in \mathcal{U} \subseteq X} \mathcal{O}_X(\mathcal{U})$$

with maps $f : X \to Y$ are pairs of a continuous map $f$ and pull-back of sections $f^{\#}_{\mathcal{U}} : \mathcal{O}_Y(\mathcal{U}) \to \mathcal{O}_X(f^{-1}\mathcal{U})$ compatible with restriction with respect to the open subsets $\mathcal{U} \subseteq Y$, and inducing **local** maps on stalk $f^{\#}_x : \mathcal{O}_{Y,f(x)} \to \mathcal{O}_{X,x}$, $f^{\#}_x(\mathfrak{m}_{Y,f(x)}) \subseteq \mathfrak{m}_{X,x}$.

The middle column consists of the category of ordinary schemes $\mathcal{S}ch$, the category of $\mathcal{Bis}$ schemes $\mathcal{BSch}$, and the category of $\mathcal{Prop}$ schemes $\mathcal{PSch}$. These are the full-subcategories of the right column consisting of $(X, \mathcal{O}_X)$ which are locally affine: we have $X = \bigcup_\alpha \mathcal{U}_\alpha$, $\mathcal{U}_\alpha$ open, and

(0.31) $$\left(\mathcal{U}_\alpha, \mathcal{O}_X\big|_{\mathcal{U}_\alpha}\right) = \mathrm{spec}\left(\mathcal{O}_X(\mathcal{U}_\alpha)\right) \qquad \text{all } \alpha.$$

These categories of generalized scheme are finite-complete and we can always construct fiber produces

(0.32)
$$\begin{array}{c} X_0 \prod_Y X_1 \\ \swarrow \quad \searrow \\ X_0 \quad \quad X_1 \\ {}_{f_0}\searrow \quad \swarrow_{f_1} \\ Y \end{array}$$

But we will be interested in more general limits of the form $\varprojlim_{j \in J} X_j$, where $J$ is a partially ordered set, and for simplicity we can restrict to

(0.33) $\quad J$ $\underline{\text{directed}}$: $\forall j_1, j_2 \in J$, $\exists j \in J$, $j \geq j_1$ & $j \geq j_2$

and

(0.34) $\quad \underline{co-finite}$ : $\forall j \in J$, $\#\{j' \in J, j' \leq j\} < \infty$,

and where $j \rightsquigarrow X_j$ is a functor, so that for $j \geq j'$ we have a map $X_j \to X_{j'}$. These limits always exist in the categories of the right column: the underlying space of $\varprojlim_{j \in J} X_j$ is the set of all coherent sequences of points $x = \{x_j\}_{j \in J}$, with the inverse limit topology, and the sheaf is the direct limit of $\mathcal{O}_{X_j}$.

Similarly, such limits always exists in the categories of the left column. Indeed we have quite formally

(0.35) $$\varprojlim_{j \in J} \operatorname{spec}(A_j) = \operatorname{spec}\left(\varinjlim_{j \in J} A_j\right)$$

But the categories of schemes, and of generalized schemes, do not have such limits: for a coherent sequence of points $x = \{x_j\}_{j \in J} \in \varprojlim_{j \in J} X_j$, while each $x_j \in X_j$ has an affine neighborhood $\mathcal{U}_j \subseteq X$, the intersection $\bigcap_{j \in J} \pi_j^{-1}(\mathcal{U}_j)$ need not be open.

We therefore pass to the pro-categories of (generalized) schemes, with objects arbitrary such inverse systems $\{X_j\}_{j \in J}$, and maps

$$(0.36) \quad \textit{pro-Sch}\,(\{X_j\}_{j \in J}, \{Y_i\}_{i \in I}) := \varprojlim_{i \in I} \varinjlim_{j \in J} \textit{Sch}\,(X_j, Y_i)$$

In the categories $\textit{pro-PSch}$ and $\textit{pro-BSch}$ we finally have the compactification of $\operatorname{spec}(\mathbb{Z})$, denoted $\overline{\operatorname{spec}(\mathbb{Z})}$, and similarly for any number field $K$ we have $\overline{\operatorname{spec}(\mathcal{O}_K)}$. For $\overline{\operatorname{spec}(\mathbb{Z})}$ we take $J$ to be the collection of finite subsets of the primes, with inclusion as partial order, and

$$(0.37) \quad X_{\{p_1 \cdots p_\ell\}} := \operatorname{spec}(\mathbb{Z}) \coprod_{\operatorname{spec}\left(\mathbb{Z}[\frac{1}{p_1 \cdots p_\ell}]\right)} \operatorname{spec}\left(\mathbb{Z}[\frac{1}{p_1 \cdots p_\ell}] \cap \mathbb{Z}_\mathbb{R}\right).$$

The generalized scheme $X_{\{p_1,\cdots,p_\ell\}}$ has closed points $\{(p_1), \cdots, (p_\ell), \eta_\mathbb{R}\}$ with $\eta_\mathbb{R} = \mathbb{Z}\left[\frac{1}{p_1 \cdots p_\ell}\right] \cap (-1, 1)$, the "**Real prime**", coming from the unique maximal ideal $(-1, 1)$ of $\mathbb{Z}_\mathbb{R}$. Enlarging the set of primes we have maps

$$(0.38) \quad X_{\{p_1,\cdots,p_\ell,q_1,\cdots,q_k\}} \longrightarrow X_{\{p_1,\cdots,p_\ell\}} \,;$$

these maps are identity on points, and are identity on the sheaves, but there are more open sets in $X_{\{p_1,\cdots,p_\ell,q_1,\cdots,q_k\}}$ then there are in $X_{\{p_1,\cdots,p_\ell\}}$. The space

$$(0.39) \quad \varprojlim \overline{\operatorname{spec}(\mathbb{Z})} \equiv \varprojlim_{\{\ell_1,\cdots,\ell_\ell\} \in J} X_{\{\ell_1,\cdots,\ell_\ell\}}$$

is one dimensional and has the generic point $(0)$, and the closed points $\eta_\mathbb{R}$ and $(p)$, where $p$ is an arbitrary prime. The structure sheaf (as a prop) is given by

$$(0.40) \quad \mathcal{O}_{\varprojlim \overline{\operatorname{spec}(\mathbb{Z})}}(\mathcal{U})_{n,m} = \left\{a \in \mathbb{Q}_{n,m}, a \circ \mathbb{Z}_p^m \subseteq \mathbb{Z}_p^n \text{ for all } p \in \mathcal{U}\right\}$$

with global sections

$$(0.41) \quad \mathcal{O}_{\varprojlim \overline{\operatorname{spec}(\mathbb{Z})}}\left(\overline{\operatorname{spec}(\mathbb{Z})}\right)_{n,m} = \mathbb{Z}_{n,m} \cap (\mathbb{Z}_\mathbb{R})_{n,m} = \mathbb{F}[\pm 1]_{n,m}$$

the $n$ by $m$ matrices where each raw and each column will have 0 entries, accept at most in one entry which could be $\pm 1$.

It is here that we obtain a truely global result: the eigenvalues of matrices in $\mathbb{F}[\pm 1]_{n,n}$ are either zero or roots of unity, the "constants" of arithmetic and conversely: the algebraic integers having "no poles over the real prime" (in any embedding into $\mathbb{C}$ they are mapped into $\mathbb{Z}_\mathbb{C}$) are the roots of unity (Kronecker Theorem).

For any prop $A$, and any $n \geqslant 1$, we have the monoid $A_{n,n}$, and its subgroup of invertible elements $\mathrm{GL}_n(A)$, the "general linear group" in dimension $n$.

For example,
(0.42)
$$\mathrm{GL}_n(\mathbb{F}) = S_n \subseteq \mathrm{GL}_n(\mathbb{F}[\pm 1]) = (\pm 1)^n \rtimes S_n \subseteq \mathrm{GL}_n(\mathbb{Z}_\mathbb{R}) = O(n) \subseteq \mathrm{GL}_n(\mathbb{Z}_\mathbb{C}) = U(n)$$

and $\mathrm{GL}_n(A)$ is the usual group of $n$ by $n$ invertible matrices for a ring $A$.

These groups come with associative maps

(0.43) $\qquad \oplus : \mathrm{GL}_{n_1}(A) \times \mathrm{GL}_{n_2}(A) \longrightarrow \mathrm{GL}_{n_1+n_2}(A)$

and give the $K$-theoretic spectrum of $A$.

The groups $\mathrm{GL}_n(A)$ act via conjugation on $A_{n,n}$, and we have the associate orbit space

$$[A_{n,n}] := A_{n,n}/\mathrm{GL}_n(A).$$

On the direct limit as $n \to \infty$ we get a commutative monoid

(0.44) $\qquad [A] = \varinjlim_{n \to \infty} [A_{n,n}] \quad , \quad [a_1] + [a_2] := [a_1 \oplus a_2]$

with an action of the multiplicative monoid of non-zero natural numbers by "Frobenius" operators

$$(0.45) \quad F_m[a] := [\underbrace{a \circ a \circ \cdots \circ a}_{m}] \quad , \quad F_{m_1} \circ F_{m_2} = F_{m_1 \cdot m_2}.$$

Applying Grothendieck's $K$-functor localizing addition, we get an Abelian group

$$(0.46) \quad \mathcal{W}(A) := K[A].$$

When $A$ is totally commutative, $A \in \mathcal{C}_T \mathcal{P}\mathit{rop}$, $\mathcal{W}(A)$ has multiplication via tensor product making it a commutative ring, and has $\lambda$-operations making it a special $\lambda$-ring.

When we take $A = \mathcal{O}_{\overline{\mathrm{spec}\mathbb{Z}}}\left(\overline{\mathrm{spec}\mathbb{Z}}\right) = \mathbb{F}[\pm 1]$, the $\lambda$-ring $\mathcal{W} = \mathcal{W}\left(\overline{\mathrm{spec}(\mathbb{Z})}\right)$ has an additive basis the cyclotomic polynomials $\phi_n$, with $\phi_1$ the unit of $\mathcal{W}$, and we have ring homomorphisms

$$(0.47) \quad t_m = \mathrm{tr} \circ F_m : \mathcal{W} = \bigoplus_{n \geq 1} \mathbb{Z} \cdot \phi_n \twoheadrightarrow \mathbb{Z}$$

with

$$(0.48) \quad t_m(\phi_n) = \mathrm{tr}(\mathbb{F}_m \phi_n) = \sum_{\xi \in \mu_n^*} \xi^m = C_n^m = \mu\left(\frac{n}{(n,m)}\right) \frac{\varphi(n)}{\varphi\left(\frac{n}{(n,m)}\right)}$$

the **Ramanujan Sums**.

The $\lambda$-ring $\mathcal{W}$ is "adic" and the operators $F_m$, $m \in \mathbb{N}^+$, extends to endomorphisms $F_m$ for

$$(0.49) \quad m \in \widehat{\mathbb{Z}}/_{\widehat{\mathbb{Z}}^*} = \prod_p p^{\mathbb{N} \cup \{\infty\}} \equiv \mathcal{CRing}(\mathcal{W}, \mathbb{Z}),$$

the super-natural-numbers, interpolating between the rank 1 projection

$F_0(\phi_n) \equiv \varphi(n) \cdot \phi_1$, and $F_1 = \mathrm{Id}_{\mathcal{W}}$:

(0.50)
$$\mathrm{tr}(\phi_n) = t_1(\phi_n) = \mu(n) \qquad \text{the Möbuis function,}$$
$$\mathrm{tr}(F_0\phi_n) = t_0(\phi_n) = \varphi(n) \qquad \text{the Euler function.}$$

We can complexify,

$$\mathcal{W}_{\mathbb{C}} = \mathbb{C} \otimes \mathcal{W} = \oplus_n \mathbb{C} \cdot \phi_n,$$

and complete $\mathcal{W}_{\mathbb{C}}$ with respect to the Hermitan form

(0.51) $\qquad < f, g >:= $ coefficient of $\phi_1$ in the product $f \cdot \overline{g}$

obtaining a Hilbert space

$$\mathcal{H} = \widehat{\mathcal{W}_{\mathbb{C}}} = \widehat{\bigoplus_n} \mathbb{C}\phi_n$$

with orthogonal basis $\phi_n$,

$$\|\phi_n\|^2 = \varphi(n).$$

The "Fourier-transform" gives an isomorphism

(0.52)
$$\mathcal{H} \xrightarrow{\sim} L_2\left(\widehat{\hat{\mathbb{Z}}}, dm\right)^{\hat{\mathbb{Z}}^*}, \quad \text{dm-additive Haar probability measure}$$

$$f \longmapsto \hat{f}(m) := \mathrm{tr}(F_m f) = t_1(F_m f) = t_m(f)$$

Considering the "zeta operators"

(0.53)
$$\zeta(F, s) = \sum_{n \geq 1} \frac{F_n}{n^s} = \prod_p \frac{1}{(1 - p^{-s} F_p)}$$

$$\zeta(F^*, t) = \sum_{m \geq 1} \frac{F_m^*}{m^t} = \zeta(F, \bar{t})^*$$

well defined on the dense subset $\mathcal{W}_{\mathbb{C}}$, and analytic for $\Re(s), \Re(t) > 1$, we have

Ramanujan's evaluations

(0.54) $\quad t_1\left(\zeta(F,s)\phi_1\right) = \zeta(s) \quad$ the Riemann's zeta function

(0.55) $\quad t_1\left(\zeta(F^*,t)\phi_1\right) = \dfrac{1}{\zeta(t)}$

(0.56) $\quad t_1\left(\zeta(F,s)\zeta(F^*,t)\phi_1\right) = \displaystyle\sum_{n,m\geqslant 1} \dfrac{C_n^m}{n^t m^s} = \dfrac{\zeta(s)}{\zeta(t)} \cdot \zeta(s+t-1).$

In chapter 1 we give the definition of $\mathcal{P}\!\mathit{rop}$s, and their commutativity. In 2 we describe the (full and faithful) embeddings of $\mathcal{R}\!\mathit{ing}$s in $\mathcal{P}\!\mathit{rop}$s. In 3 we give the definition of $\mathcal{B}\!\mathit{io}$s, chapter 4 discuss their commutativity. In Chapter 5 we describe the ideals and primes of a commutative $\mathcal{B}\!\mathit{io}$, $A$, and in chapter 6 we describe the (compact, sober) Zariski topology on spec($A$). In chapter 7 we give the basic facts for localizations, and in chapter 8 construct the sheaf of $\mathcal{B}\!\mathit{io}/\mathcal{P}\!\mathit{rop}$s over spec($A$). In chapter 9 we describe the categories of $\mathcal{B}\!\mathit{io}/\mathcal{P}\!\mathit{rop}$s schemes, and in chapter 10 we describe the pro-categories of these, and the "compactification" $\overline{\mathrm{spec}\mathcal{O}_K}$, $K$ a number field. In chapter 11 we define valuation props and bios, we show that they correspond in the rank 1 case to real valued valuations; for a number field $K$ they correspond to the finite primes $\mathfrak{p} \subseteq \mathcal{O}_K$ or to the real/complex primes $\sigma : K \hookrightarrow \mathbb{C}$, $\sigma \sim \overline{\sigma}$. We also show how the operad structure is compatible with Beta-integrals, and give the "zeta machine": given a compact homogeneous valuation $\mathcal{P}\!\mathit{rop}$ $B$, together with a map $\mathbb{N} \to K = (B(1))\backslash\{0\})^{-1}(B) \equiv$ fraction field of $B$, we get a meromophic function $L(B,s)$, normalized by $L(B,1) = 1$ ; we have $L(\mathbb{Z}_p,s) = \dfrac{\zeta_p(s)}{\zeta_p(1)}$, with $\zeta_p(s) = (1-p^{-s})^{-1}$, $\zeta_\mathbb{R}(s) = 2^{\frac{s}{2}}\Gamma(\frac{s}{2})$, $\zeta_\mathbb{C}(s) = \Gamma(s)$. In chapter 12 we describe for a $\mathcal{B}\!\mathit{io}/\mathcal{P}\!\mathit{rop}$ $A$, and a category $\mathcal{C}$, the $A$-objects in $\mathcal{C}$. In chapter 13 we sin again with addition, taking $\mathcal{C} = \mathcal{A}\!\mathit{b}$ we have the category of $A$-modules, isomorphic to a category of abelian group objects in the category of $\mathcal{Bio}/\mathcal{Prop}$ over $A$ giving rise to derivations, and to the Kahler differentials representing them. These all arise from our attempt to linearize geometry infinitesimally. We describe explicitly the $\mathbb{Z}$-prop-module $\Omega(\mathbb{Z}/\mathbb{F})$. In chapter 14 we bring in the simplicial $\mathcal{Bio}/\mathcal{Prop}$, and the Quillen model structure on them, giving rise to the cotanget complex of a homomorphism of $\mathcal{Bio}/\mathcal{Prop}$. In chapter 15 we define briefly the basic properties of maps of $\mathcal{Bio}/\mathcal{Prop}$ schemes. In chapter 16 we give for a $\mathcal{Prop}$ $A$, commutative or not, the infinite loop space $\mathcal{K}_A$ whose homotopy groups are the higher $K$-groups of $A$; we do it by a very concrete version of Quillen's $S^{-1}S$ construction. In chapter 17 we give the Witt ring story, and the basic example of the special $\lambda$-ring $\mathcal{W} = \mathcal{W}(\overline{\mathrm{spec}\mathbb{Z}})$. In chapter 18 we give the (close) symmetric monidal structure on (totally-commutative) $A$-sets for a $\mathcal{Bio}$ $A$, and we briefly describe the stabilization of simplicial $A$-sets using "symmetric spectra" or modules over the sphere spectrum - the point being is (Jeff Smith idea [**HSS00**]) that the spheres $S^{\cdot} = \{S^n\}_{n \geq 0} \in (s\mathcal{Set})^{\amalg_n S_n} \equiv \sum(\mathbb{F})$ are a **commutative** monoid - this shows very clearly the interaction of homotopy and arithmetic via $\mathbb{F} \subseteq \mathbb{Z}$.

We try to avoid most technicalities, proofs or calculations, and we concentrate on the concepts and ideas. Two exceptions are the calculation of the global sections of the sheaf associated with a prop or a bio (Theorem 8.1), and the proof of Ostrowski theorem- the calculation of the valuations on a number field (Theorem 11.2). We included these so that the reader can see that they are exactly the same as the proofs of the analogue classical statements only formulated in the language of props or bios (we suggest the reader skip these proofs on a first reading). For the same reason we avoid the globalizations of chapters (16) and (18), which are best done in the more technical language of infinite categories.

This monograph is dedicated to the memory of my teachers, mentors, and friends

- Daniel Quillen
- Yuri Manin
- John Coates